\UseRawInputEncoding
\documentclass[12pt,twoside]{amsart}
\usepackage{amsmath}
\usepackage{epsfig}
\usepackage{amssymb}
\usepackage{xcolor}
\usepackage{blkarray}

\usepackage{tikz}
\usetikzlibrary{calc}
\usetikzlibrary{intersections}
\usepackage{tikz-3dplot}

\newtheorem{theorem}{Theorem}[section]
\newtheorem{lemma}[theorem]{Lemma}

\newtheorem{claim}{}[theorem]

\newtheorem{conjecture}[theorem]{Conjecture}

\theoremstyle{definition}

\theoremstyle{remark}

\numberwithin{equation}{section}

\newcommand{\si}{{\rm si}}
\newcommand{\cl}{{\rm cl}}

\newcommand{\cB}{{\mathcal B}}

\newcommand{\subproof}{\begin{proof}[Subproof]}

\newcommand{\del}{\backslash}
\newcommand{\con}{/}
\newcommand{\ret}{\rightharpoonup}
\newcommand{\proj}{/\!\! /}

\begin{document}

\sloppy

\title{Triples in matroid circuits}

\author[Geelen]{Jim Geelen}
\thanks{Jim Geelen (Department of Combinatorics and Optimization,
University of Waterloo,  Canada) was
supported by grants from NSERC and the Office of Naval Research. }

\author[Gerards]{Bert Gerards}

\thanks{Bert Gerards (Sittard, The Netherlands).}

\author[Whittle]{Geoff Whittle}

\thanks{Geoff Whittle (School of Mathematics and Statistics, Victoria University of Wellington, New Zealand)
was supported by a grant from the Marsden Fund of New Zealand.}

\thanks{This project is also supported through a generous personal donation from Neil Robertson to the University of Waterloo.}


\begin{abstract}
We prove that, if $X$ is a four-element independent set in a vertically $5$-connected representable matroid $M$, then either there is a circuit of $M$ that contains at least three elements of $X$, or there is a representation $A$ of $M$ such that each column of $A$ has at most two non-zero elements and each of the columns indexed by $X$ has exactly one non-zero element. We also prove a version of this result without the connectivity assumptions. 
\end{abstract}

\maketitle

\section{Introduction}

A {\em frame matrix} is a matrix in which each column has at most two non-zero entries. A {\em unit-column} is a column with exactly one non-zero element.
If a matroid $M$ is represented by a frame matrix $A$ and $X$ is the set of elements in $M$ that index unit columns, then no circuit in $M$ can contain three or more elements in $X$. We prove the following partial converse. Here two matrices are {\em row-equivalent} if they have the same row-space.

\begin{theorem}
\label{5conn}
Let $M$ be a vertically $5$-connected matroid represented by a matrix $A$. If $X$ is an independent set, 
with $|X|\ge 4$, then either there is a circuit of $M$ that contains at least three elements of $X$, or $A$ is row-equivalent to a frame matrix $A'$ and each element of $X$ indexes a unit column of $A'$.
\end{theorem}

Vertical $5$-connectivity is necessary, which we discuss in Section~2, but our main result has no connectivity restrictions. Roughly, Theorem~\ref{one} states that, if $X$ is an independent set, with $|X|\ge 4$, in a representable matroid  $M$  and there is no circuit containing three or more elements of $X$, then $M$ can be constructed from a frame matroid (as in Theorem~\ref{5conn}) by attaching other matroids via $1$-, $2$-, $3$-, and $4$-sums. Moreover, the $4$-sums attach onto planes spanned by three unit columns. This result is an important tool in our matroid minors project, which we discuss in the following section.

\section{Motivation}\label{sec:motivation}

We call a set $X$ of elements in a matroid $M$ {\em coupled} if there is no circuit in $M$ that contains three of more elements of $M$.
Note that $X$ is coupled if and only if there is no $U_{2,3}$-minor $N$ of $M$ with $E(N)\subseteq X$. Therefore, if $X$ is coupled in $M$ then it is also coupled
in any minor of $M$ that contains $X$.

Given three distinct elements $a$, $b$, and $c$ in a matroid $M$,
the {\em three-elements problem} is the problem of determining whether or not $\{a,b,c\}$ is coupled.

\subsection*{Binary matroids}
The three-elements problem was solved for cographic matroids by
Chakravarti and Robertson~\cite{Neil} and then extended to binary matroids
by Seymour~\cite{Seymour:triples} who proved the following result.
\begin{theorem}\label{binarytriples}
If $\{a,b,c\}$ is a coupled set in a vertically $4$-connected, binary matroid $M$, then $M$ is the cycle-matroid of a graph $G$ and $a$, $b$, and $c$ are incident with a common vertex in $G$.
\end{theorem}
Note that if $v$ is the vertex of $G$ that is incident with $a$, $b$, and $c$ and $A$ is the matrix obtained from 
a  signed incidence matrix of $G$ representing $M$ by deleting the row indexed by $v$, then $A$ is a frame matrix and each of $a$, $b$, and $c$ indexes a unit column of $A$.

Theorem~\ref{binarytriples} provides an efficient algorithm for solving the three-elements problem in a binary matroid $M=M(A)$.
Note that, adding points to the lines spanned by $\{a,b\}$, $\{b,c\}$, and $\{a,c\}$ does not affect whether or not $\{a,b,c\}$ is coupled; therefore we may assume that $M$ has an $M(K_4)$-restriction $N$ containing $a$, $b$, and $c$. In the case that $M$ is vertically $4$-connected, we can test whether or not $\{a,b,c\}$ is coupled by simply testing whether $M$ is graphic. Indeed, if $M$ is graphic, then the minor $N$ will force $a$, $b$, and $c$ to be incident with a common vertex. In the case that $M$ is not vertically $4$-connected, there is a vertical $3$-separation $(X,Y)$ with $a,b,c\in X$. Now we can construct a binary extension of $M$ by adding (all possible) points spanned by both $X$ and $Y$, then deleting $Y$. It is routine to show that this does not affect whether or not $\{a,b,c\}$ is coupled (assuming that we decompose on vertical separations of minimum order first). This decomposition decreases the rank, so repeating in this way we will eventually reduce to the vertical $4$-connected case.

Our main result (Theorem~\ref{one}) is in fact a generalization of Theorem~\ref{binarytriples}. To see the connection, we can construct a binary matroid $M^+$ by coextending $M$ by a new element $d$ so that $\{a,b,c,d\}$ is a cocircuit in $M^+$. In binary matroids any circuit meets any cocircuit in an even number of elements. Therefore, $\{a,b,c\}$ is coupled in $M$ if and only if $\{a,b,c,d\}$ is coupled in $M^+$. It still requires some additional work to obtain Theorem~\ref{binarytriples}, but the key observation is that, since $\{a,b,c,d\}$
is a cocircuit, the decomposition given by Theorem~\ref{one} will not involve $4$-sums.

\subsection*{Intractability for general matroids}
If $T$ is a triangle in a $3$-connected  non-binary matroid, then there is a $U_{2,4}$-minor containing $\{a,b,c\}$; see Theorem~\ref{modular triangle}. Note that $T$ will be a cocircuit in such a minor. Dually this gives the following result.
\begin{theorem}\label{thm:triad}
If $\{a,b,c\}$ is a triad in a $3$-connected matroid $M$, then $\{a,b,c\}$ is coupled in $M$ if and only if $M$ is binary.
\end{theorem}
This result is interesting from a structural point-of-view, but it has some shortcomings with respect to algorithms.
Seymour~\cite{Seymour:nonbinary} showed that the problem of deciding whether or not a matroid is binary is intractable, under the rank-oracle model, and we will use a variation of that argument to show that the three-elements problem is also intractable.

For any integer $n\ge 3$, we take $M_n$ to be a rank-$(n+1)$ ``binary spike with tip  and co-tip".
Specifically,  $M_n$ is a binary matroid with ground-set $\{b_0,b_1,\ldots,b_n\}\cup\{a_0,a_1,\ldots,a_n\}$ where
 $\{b_0,b_1,\ldots,b_n\}$ is a basis, each of $\{b_0,b_1,a_1\},\ldots,\{b_0,b_n,a_n\}$ is a triangle, and
each of $\{a_0,b_1,a_1\},\ldots,\{a_0,b_n,a_n\}$ is a triad.
It is easy to check that, $M_n$ has at least $2^{n-1}$ circuit-hyperplanes and that relaxing any one of these circuit-hyperplanes will result in a non-binary matroid. Now $\{a_0,a_1,b_1\}$ is not contained in a circuit in $M_n$, but it is in each of the relaxations. Therefore we cannot safely conclude that $\{a_0,a_1,b_1\}$ is coupled without explicitly checking the rank of exponentially many sets. That is, the three-elements problem cannot be solved in polynomial-time under the rank-oracle model.   

We expect that the three-elements problem can be solved in polynomial-time when $M$ is given by a matrix representation. We also hope to find a complete structural characterization for coupled triples in represented matroids.

\subsection*{Vertical $5$-connectivity is necessary} Consider the matroid $M$ represented by the following matrix over any field with characteristic different from $2$,
$$
\begin{array}{c}
\begin{array}{ccc}
\hspace{2cm}\color{blue} B\hspace{1cm} & \hspace{1cm}\color{red} X \hspace{1cm}&\hspace{2cm} Y \hspace{4cm}
\end{array}\\
\left[
\begin{array}{rrrr|rrr|rrrrrr}
\color{blue} 1 & \color{blue} 0 & \color{blue} 0 & \color{blue}  0&
\color{red} 1  & \color{red} 1  & \color{red} 1  & 1 & 1 &  1 & 0 & 0 & 0  \\
\color{blue} 0 & \color{blue} 1 & \color{blue} 0 & \color{blue} 0 &
\color{red} 1  & \color{red} -1  & \color{red} -1  & -1 & 0 &  0 & 1 & 1 & 0  \\
\color{blue} 0 & \color{blue} 0 & \color{blue} 1 & \color{blue} 0 &
\color{red} -1  & \color{red} -1  & \color{red} 1  & 0 & -1 &  0 & -1 & 0 & 1  \\
\color{blue} 0 & \color{blue} 0 & \color{blue} 0 & \color{blue} 1 &
\color{red} -1  & \color{red} 1  & \color{red} -1  & 0 & 0 &  -1 & 0 & -1 & -1  
\end{array}\right].
\end{array}
$$
Note that
\begin{itemize}
\item[(a)] $M|(B\cup Y)$ is isomorphic to $M(K_5)$,
\item[(b)] $B$ is a cocircuit in $M$ and every pair of elements in $B$ spans an element in $Y$,
\item[(c)] every element in $Y$ is spanned by two elements in $X$, and 
\item[(d)] $M|(B\cup X)$ is the dual of the non-Fano.
\end{itemize}
Using $(c)$, the set $X$ is coupled in $M|(X\cup Y)$ and, by $(b)$, it remains coupled in $M$.
Moreover, $M$ is vertically $5$-connected, since $M$ is spanned by $M(K_5)$.
We leave it to the reader to check that $M$ cannot be represented by a frame matrix.
We can use $4$-sums to attach these matroids to frame representations to show that the $5$-connectivity requirement in Theorem~\ref{5conn} is necessary. 

The problem with this example is that there is a frame representation of $M(K_4)$ that 
does not extend to a frame representation of $M(K_5)$, so we can take a $4$-sum of a frame representation and $M(K_5)$ such that the resulting matroid is not a frame representation. 
More generally, we can do the same thing for $M(K_n)$, when $n\ge 5$. There are other constructions, but a complete characterization, in the representable case,  still looks achievable.

\subsection*{The two-linkage problem} 
Given distinct vertices $s_1, \,t_1,\, s_2,\, t_2$ in a graph $G$, the {\em two-linkage problem} is the problem of deciding whether or not there exist two disjoint paths, one connecting $s_1$ to $t_1$ and the other connecting $s_2$ to $t_2$. This is in fact a special case of the three-elements problem for binary matroids. We may add edges $a=s_1 s_2$, $b=s_1t_2$, and $c=s_2t_2$ without affecting whether or not the linkage exists. It is an easy exercise to verify that the linkage exists if and only if the elements $a$, $b$, and $c$ are contained in a  circuit in $M(G)^*$ . Using this, one can see that Seymour's theorem (\ref{binarytriples}) implies the following result that was proved independently by
Seymour\cite{Seymour:linkage} and Thomassen~\cite{Thomassen}.
\begin{theorem}\label{two-linkage}
Let $s_1, \,t_1,\, s_2,\, t_2$ be disjoint vertices in a $4$-connected graph $G$.
If there do not exist two disjoint paths one connecting $s_1$ to $t_1$ and the other connecting $s_2$ to $t_2$, then $G$ admits a planar
embedding with $s_1, \,s_2,\, t_1, \, t_2$ embedded in that order on the outer face.
\end{theorem}

\subsection*{Graph minors}
The two-linkage problem plays an important role in the proof of the graph-minors structure theorem~\cite{GMXVI} of Robertson and Seymour. In particular, Theorem~\ref{two-linkage} is used in the proof of the ``flat-wall theorem''~\cite{GMXIII} which roughly states that, for $k_2\gg k_1\gg t\gg n$, given a $k_2\times k_2$-grid-minor in a graph $G$ with no $K_n$-minor, there is a subgraph $H$ of $G$ such that $H$ embeds in the plane ``up to $3$-separations", $H$ includes a $k_1\times k_1$-sub-grid of our grid minor, and there is a set $X\subseteq V(G)\setminus V(H)$ of at most $t$ vertices such that no vertex of $G$ outside $H$ and $X$ is connected to any vertex of $H$ except those on the outer face in the planar embedding. 

Our main result plays a similar role in structure theory for matroids represented over a finite field.

\section{A Sylvester-Gallai Theorem for planes}\label{section SG}
As usual, we use the terms {\em point}, {\em line}, and {\em plane} for a flat of rank one, two, or three respectively.
In particular, note that a point is a set, not an element.
We say that two lines in a matroid are {\em skew} if their union has rank $4$, and we call a line {\em long} if it has at least three points.

The Sylvester-Gallai Theorem states that every rank-$3$ real-representable matroid has a two-point line; see~\cite{Sylvester, Gallai, Melchior}.
Motzkin~\cite{Motzkin} observed that this result does not directly extend to rank-$4$, since the union of two skew long lines has no three-point plane.
The problem of completely characterizing the rank-$4$ real-representable matroids that have no three-point plane remains open,
but Motzkin proved a significant partial result, stated below.

A plane is called {\em ordinary} if it is the union of a line and a point. Thus, in a matroid without any long lines, an ordinary plane has only three points.
\begin{lemma}[Motzkin's Lemma]\label{Motzkins lemma}
Every rank-$4$ real-representable matroid contains an ordinary plane.
\end{lemma}

In Section~\ref{frameworks}, we will apply Motzkin's Lemma to gammoids, which are matroids that arise from connectivity in directed graphs.
Specifically, a matroid $M$ on a set $Z$ is a {\em gammoid} if there is a 
directed graph $G=(V,E)$, with $Z\subseteq V$, and a set $X\subseteq V$ such that the rank of  each $W\subseteq Z$ is equal to the 
maximum number of disjoint directed paths from $W$ to $X$ in $G$. The class of gammoids is minor-closed; in fact a matroid is a gammoid if and only if it is a minor of a transversal matroid; see Oxley~\cite[Section 3.2]{Oxley}.
Since transversal matroids are real-representable~\cite[Corollary 11.2.18]{Oxley}, gammoids are also real-representable. 

We guess that gammoids satisfy the following strengthening of Motzkin's Lemma.
\begin{conjecture}
Every rank-$4$ gammoid either has a three-point plane or it is the union of two skew long-lines.
\end{conjecture}

\section{Connectivity in undirected graphs}\label{graph connectivity}

The gammoids that we are interested in come from connectivity in {\em undirected} graphs; we can replace each edge 
with a pair of oppositely directed edges without affecting the connectivities between any two sets of vertices, so these are, in fact, gammoids.

We take the convention that the empty graph is connected. In a departure from matroid terminology, $3$-connected graphs are not required to be simple.

For subgraphs $A$ and $B$ of a graph $G$, we can define $A\cup B$ and $A\cap B$, naturally, by taking intersections or unions of both the vertex sets and edge sets.
A {\em separation} in a graph $G$ is a pair $(A,B)$ of subgraphs of $G$ such that $G=A\cup B$; the {\em order} of the separation is $|V(A\cap B)|$.  We call a separation $(A,B)$ {\em proper} if $V(A)\setminus V(B)$ and $V(B)\setminus V(A)$ are both non-empty.
For sets $X$ and $Y$ of vertices, an {\em $(X,Y)$-separation} is a separation $(A,B)$ with $X\subseteq V(A)$ and $Y\subseteq V(B)$.
We define $\kappa_G(X,Y)$ to be the minimum order of an  $(X,Y)$-separation. Note that we do not require $X$ and $Y$ to be disjoint.

Let $X$ be a set of vertices in a graph $G=(V,E)$. As discussed above, there is a gammoid $K(G,X)$ on $V$ such that the rank of a set $W$ is $\kappa_G(W,X)$. 

\subsection*{Tight sets}
Let $e$ be an edge in a graph $G$. We write $G_{e\ret z}$ to denote the graph obtained from $G$ by contracting $e$ to a vertex $z$. 
Let $X$ and $Z$ be sets of vertices in a graph $G=(V,E)$. We say that $Z$ is {\em $X$-tight}, if, for each vertex $v\in V\setminus Z$, the 
connectivity from $Z\cup\{v\}$ to $X$ is $|Z|+1$. In terms of the connectivity matroid, $K(G,X)$, this says that $Z$ is both independent set and a flat. In particular, an $X$-tight triple is a three-point plane in $K(G,X)$.

The notion of $X$-tightness is very useful in matroid minors. For example, one can show that, if $(G,T)$ is a graft and $Z$ is an $X$-tight subset of $T$, then we can contract all of $G$ down to $X$ keeping the $T$-set of size at least $|Z|$. That is essentially 
the content of Lemma~\ref{extend framework 2}, except that the lemma is not restricted to binary matroids.

The following result shows how to remove  edges that are not incident with $Z$. There is a companion result, below, for removing  edges incident with $Z$. 
\begin{lemma}\label{tight sets 1}
Let $X$ and $Z$ be sets of vertices in a graph $G=(V,E)$, let $e=xy$ be an edge of $G-Z$ with  $y\not\in X$.
If $Z$ is $X$-tight in $G$, then $Z$ is $X$-tight in at least one of $G-e$ or $G_{e\ret x}$. 
\end{lemma}

\begin{proof}
Let $k=|X|$. If $Z$ is not $X$-tight in $G_{e\ret x}$, then there is a $(Z,X)$-separation $(A_1,B_1)$ of order $k+1$ in $G-e$  such that
$\{x,y\}\subseteq V(A_1\cap B_1)$. If $Z$ is not $X$-tight in $G\del e$, then there is a $(Z,X)$-separation $(A_2,B_2)$ of order $k$ in $G-e$  such that $y\in V(A_2)\setminus V(B_2)$ and $x\in V(B_2)\setminus V(A_2)$. Note that
$$ (A_1\cap A_2, B_1\cup B_2) \mbox{ and } (A_1\cup A_2, B_1\cap B_2)$$
are both $(Z,X)$-separations in $G-e$ and, by submodularity, the sum of the orders of these separations is at most $2k+1$.
Since $\{x,y\}\subseteq V(A_1\cap B_1)$,
$$ (A_1\cap A_2, (B_1\cup B_2)+e) \mbox{ and } ((A_1\cup A_2)+e, B_1\cap B_2)$$
are both $(Z,X)$-separations in $G$ and the sum of the orders of these separations is at most $2k+1$.
However, since $x$ is a vertex  in both $A_1\cap A_2$ and $(A_1\cup A_2)+e$, this contradicts the fact that 
$Z$ is $X$-tight in $G$.
\end{proof}

Next we show we can reduce the degree of each vertex in $Z\setminus X$ to two. 
\begin{lemma}\label{tight sets 2}
Let $X$ and $Z$ be sets of vertices in a graph $G=(V,E)$, and let $z\in Z\setminus X$. If $Z$ is $X$-tight and $z$ has degree at least $3$, then there is an edge $e$ incident with 
$z$ such that $Z$ is $X$-tight in $G-e$.
\end{lemma}

\begin{proof} Let $k=|X|$.
Consider constructing a new graph $G^+$ by adding a twin $z'$ of $z$ to $G$. Since  $Z$ is $X$-tight in $G$, there exist
$k+1$ disjoint $(Z\cup\{z'\},X)$-paths in $G^+$. Therefore there are $(Z,X)$-paths $P_0,P_1,\ldots,P_k$ such that
$P_0\cup P_1,P_2,\ldots,P_k$ are disjoint and $V(P_0\cap P_1)=\{z\}$. Suppose that $e=zy$ is an edge incident with $z$ that is in neither 
$P_0$ nor $P_1$. We claim that $Z$ is $X$-tight in $G-e$. Suppose otherwise and consider an $(Z,X)$-separation $(A,B)$ of order $k$ in $G-e$
with $V(A)\neq Z$. Since $Z$ is $X$-tight in $G$, neither $A$ nor $B$ can contain both ends of $e$. Thus
$z\in V(A)\setminus V(B)$ and $y\in V(B)\setminus V(A)$. This is not possible, since the $k+1$ paths $P_0,P_1,P_2,\ldots,P_k$ must each intersect $V(A\cap B)$ at distinct vertices.
\end{proof}

\section{Matroid connectivity}\label{matroid connectivity}

We assume that the reader is familiar with matroid theory; we mostly follow the terminology and notation of Oxley~\cite{Oxley}. There is a significant interplay between graphs and matroids in this paper. To make the connection relatively seamless, we use non-standard terminology for ``separations'' in a matroid.

Recall that the {\em local connectivity} between two sets $X$ and $Y$ in a matroid $M$ is defined as 
$$ \sqcap_M(X,Y) := r(X) + r(Y) - r(X\cup Y).$$
The {\em Tutte connectivity} of a set $X$ is
$$ \lambda_M(X) := \sqcap_M(X,E(M)\setminus X).$$
A {\em separation} is a pair $(A,B)$ of sets in $M$ such that
$A\cup B=E(M)$; the {\em order} of the separation is $\sqcap_M(A,B)$.
For two sets $X$ and $Y$  in a matroid $M$, an {\em $(X,Y)$-separation} is a separation $(A,B)$ 
with $X\subseteq A$ and $Y\subseteq B$.
The {\em connectivity between} two sets $X$ and $Y$, denoted by $\kappa_M(X,Y)$, is defined as the minimum order of an $(X,Y)$-separation in $M$.
Note that this differs from the definition on Oxley~\cite[page 320]{Oxley} since we do not require that $X$ and $Y$ are disjoint. Our definition is tantamount to making $X$ and $Y$ disjoint by adding parallel copies of the common elements
and then applying the definition in Oxley. 

We recall that a {\em vertical separation} in a matroid $M$ is a {separation} 
$(A,B)$ of $E(M)$ such that neither $A$ nor $B$ is spanning.
We say that $M$ is {\em vertically $k$-connected} if it has no vertical separation of order less than $k-1$.

Note that, Oxley's ``$k$-separations" are separations of order at most $k-1$ here.
Nevertheless, our usage of ``$k$-connected" and ``vertically $k$-connected" are consistent with Oxley's definition.

For ease of reference, here are two well-known connectivity results that we use. 
\begin{theorem}[{Oxley~\cite[Theorem 4.3.1]{Oxley}.}]\label{delcon}
If $e$ is an element in a connected matroid $M$, then $M\del e$ or $M\con e$ is connected. 
\end{theorem}

\begin{theorem}[{Tutte's linking theorem; see Oxley~\cite[Theorem 8.5.2]{Oxley}.}]
 If $X$ and $Y$ are disjoint sets of elements in a matroid $M$ and $e\in E(M)\setminus (X\cup Y)$, then 
 $$\max(\kappa_{M\del e}(X,Y),\, \kappa_{M\con e}(X,Y))=\kappa_M(X,Y).$$
 \end{theorem}

\subsection*{Strands}
Let $N$ be a restriction of a matroid $M$. We call a set $S\subseteq E(M)$ a {\em strand}
for $N$ if  $\sqcap_M(E(N),S)=1$ and each pair of elements in $S$ is in series in 
 $M|(E(N)\cup S)$. The following result will be familiar to anyone who has thought about ear decomposition for
 matroids; Tutte~\cite{Tutte:lectures} uses similar results in his algorithm for testing whether a binary matroid is graphic.

\begin{lemma}\label{strand}
If $N$ is a restriction of a matroid $M$ and $e\in E(M)\setminus E(N)$ is such that the
component of $M$ that contains $e$ also contains at least one element of $N$,
then there is a strand for $N$ in $M$ that contains $e$.
\end{lemma}

\begin{proof}
Suppose otherwise, and consider a counterexample with $|M|$ minimum.
Let $C$ be a circuit of $M$ that contains $e$ and at least one element of $N$.
We may assume that $e$ is not spanned by $N$, since otherwise $\{e\}$ is already a strand.
Thus there is some element, say $f\in C\setminus\{e\}$ that is not spanned by $N$.
By minimality, there is a strand $S$ for $N$ in $M\con f$ that contains $e$.
But then, either $S$ or $S\cup \{f\}$ is a strand for $N$ in $M$.
\end{proof}

\section{Modular flats}\label{modular flats}

A pair $(F_1,F_2)$ of subsets of $E(M)$ is {\em modular} if  $r(F_1\cap F_2) =\sqcap_M(F_1,F_2)$.
A  set $F$ of elements is {\em modular} if $(F,F')$ is a modular pair for each flat $F'$ of $M$. In a simple matroid a modular set is necessarily a flat.
It is easy to see that modularity is preserved under minors.
That is, if $F$ is a modular flat and $e\in E(M)\setminus F$, then $F$ is modular in both $M\del e$ and $M\con e$. Moreover, if $e\in F$, then $F\del\{e\}$ is also modular in $M\con e$, but deleting elements of $F$ need not preserve modularity. The following result is well-known; see, for example Kapadia~\cite[Proposition 1.4.1]{Kapadia}.
\begin{lemma}\label{modular sets}
For a set $F$ of elements in a matroid $M$, the following are equivalent:
\begin{enumerate}
\item[(1)]
$F$  is a modular set in $M$,
\item[(2)]
for each minor $N$ of $M$ containing $M|F$ as a spanning restriction, each element of $E(N)\setminus F$
is either a loop or is in parallel with an element of $F$.
\item[(3)] each strand of $M|F$ in $M$ spans a non-loop element of $F$.
\end{enumerate}
\end{lemma}

\begin{proof}
If $(1)$ fails, then  there is a flat $F'$ of $M$ such that $(F,F')$ is not a modular pair. Let $I$  be a maximal independent set of $(M\con F)|F'$and let $N=M\con I$. Note that $F'\setminus I$ is contained in the closure of $F$ in $N$ and $r_N(F\cap F')=r_M(F\cap F')< \sqcap_{M}(F,F') = r_N(F'\setminus I)$ so there is an element $e\in F'\setminus I$ that is not in parallel with any element of $F$, so $(2)$ fails.

If $(2)$ fails, there is a maximal independent set $I$ of $M\con F$ and an element $e\in E(M)\setminus I$, such that, in $M\con I$, $e$ is not a loop and is not in parallel with any element of $F$. Let $Y:=I\cup\{e\}$.
Then $\sqcap_M(F,Y)=1$ and yet $Y$ does not span a non-loop element of $F$.
Since $\sqcap_M(F,Y)=1$, there is a circuit of $M|(F\cup Y)$ that 
contains an element of $F$ and an element of $Y$. By Lemma~\ref{strand},
$Y$ contains a strand for $M|F$ and, hence, $(3)$ fails.

If $(3)$ fails and $S$ is a strand of $M|F$ in $M$ that does not span any non-loop element of $F$,
then $(F,\cl_M(S))$ is not modular, so $(1)$ fails, as required.
\end{proof}

A matroid is {\em round} if its ground-set cannot be covered by two hyperplanes. The next result was probably already known to Tutte~\cite{Tutte:lectures} who used similar results in his algorithm for testing whether a binary matroid is graphic.
\begin{lemma}\label{modular implies round}
Let $X$ be a modular set of elements in a matroid $M$ and let $N$ be a component of $M\con X$. If $Z$ is the set of elements in $X$ that are spanned by $E(N)$ in $M$, then $M|Z$ is round.
\end{lemma}

\begin{proof} Suppose that the result fails, and consider a counterexample $(M,X,N)$ with $|M|$ minimum.
By minimality, we have $E(M)=X\cup E(N)$, $X=Z$, and $M$ is simple and connected.
Since $(M,X,N)$ is a counterexample,  there exist hyperplanes $H_1$ and $H_2$ of $M|X$ such that $X=H_1\cup H_2$.
Note that  $H_1$ and $H_2$ are disjoint since otherwise we could contract the intersection. 
Consider an element $e$ in $N$.
By minimality $N\con e$ is not connected, so, by Theorem~\ref{delcon}, $N\del e$ is connected.
Again applying minimality, $E(N)\del e$ cannot span an element in both $H_1$ and $H_2$, so, up to symmetry,
 $\sqcap_M(E(N\del e), H_1)=0$. But then $\sqcap_M(E(N), H_1)=1$ and hence $|H_1|=1$. Since $H_1$ and $H_2$ are both hyperplanes of
 $M|X$, we also have $|H_2|=1$. Since $M$ is connected there is a circuit $C$ that contains both elements of $X$. Now $C\setminus X$
 is a strand for $M|X$, but $C\setminus X$ spans no  element of $X$, contrary to the fact that $X$ is modular.
\end{proof}

\subsection*{Structure theory}
Modular flats have an interesting connection with structure theory. The first result is an easy consequence of Seymour's result~\cite{Seymour:nonbinary} that $U_{2,4}$ is $2$-rounded; see Oxley~\cite[Proposition 12.3.8]{Oxley}.
The result was stated in this form in Rohan Kapadia's PhD thesis~\cite[Theorem 1.4.2]{Kapadia}.
\begin{theorem}\label{modular triangle}
If a triangle $T$  in a $3$-connected matroid $M$ is modular, then $M$ is binary.
\end{theorem}

\begin{proof} Let $T=\{a,b,c\}$. If $M$ is not binary, $M$ has a $U_{2,4}$-minor. Then, by the $2$-roundedness of $U_{2,4}$,
there is a $U_{2,4}$-minor $N$ that contains $\{a,b\}$. We may assume that $N$ also contains $c$; to see this  consider what happens when $c$ is the last element removed when constructing $N$. But then, by Lemma~\ref{modular sets}, $T$ is not modular in $M$.
\end{proof}

Here is an even more striking example, which is just Seymour's Theorem~\ref{binarytriples} in disguise; this result first appeared in this form in~\cite[Theorem 1.4.3]{Kapadia}.
\begin{theorem}\label{modular K4}
Let $N\cong M(K_4)$ be a  restriction of a vertically $4$-connected matroid $M$. If $E(N)$ is modular in $M$, then $M$ is graphic.
\end{theorem}

\begin{proof}  Each triangle in $N$ is modular in $M$, so, by Theorem~\ref{modular triangle}, $M$ is binary. Now we obtain the result by applying Theorem~\ref{binarytriples} to a triad $\{a,b,c\}$ of $N$.
\end{proof}

\subsection*{Modular sums} The following result is easy and well-known; we omit the proof.
Here, we are primarily interested in the case that $N=U_{2,3}$.
\begin{lemma}\label{decomp} 
Let $N$ be a matroid and let $(F_1,F_2)$ be a modular pair of flats in a matroid $M$ with $E(M)=F_1\cup F_2$.
If $E(N)\subseteq F_1$, then $N$ is a minor of $M$ if and only if $N$ is a minor of $M|F_1$.
\end{lemma}
The matroid $M$  is a generalized parallel connection (or {\em modular sum}) of $M|F_1$ and $M|F_2$; see Oxley~\cite[Section 11.4]{Oxley}. The lemma 
shows that we can use modular sums to build up matroids that avoid $N$ as minor.
We assume that $M\con F_1$ is connected, since otherwise we could construct $M$ via a series of modular sums, treating the components one at a time.  

Note that, by Lemma~\ref{modular implies round}, the matroid $M|(F_1\cap F_2)$ is round. The cycle matroid of a graph $G$ is round if and only if $G$ is complete. This is why we only see clique-sums in the graph minors structure theory. Roundness is significantly less prescriptive for general matroids. For example, all $3$-connected rank-$3$ matroids are round, giving rise to the possibility of other interesting $4$-sums.

\section{Planted matroids and modularity}

One useful property of representable matroids is that when $\sqcap_M(X,Y)>0$ we can extend $M$ by a non-loop element $e$ in the span of both $X$ and $Y$.
We rely more on such extension properties  than on specific coordinatizations. 
For this reason it is convenient to think of a representation 
of a simple matroid as a restriction of a (possibly infinite) projective geometry of finite rank. Moreover, we will also consider representations of non-simple matroids, so we need to add loops and parallel elements to the projective geometry.    

We are almost exclusively interested in finite matroids, but for the following brief discussion we allow infinite matroids of finite rank.
We say that a matroid is {\em modular} if each pair $(F_1,F_2)$ of flats is modular.
Projective geometries are modular and any matroid that simplifies to a projective geometry is also modular.
A {\em modular space} is a (possibly infinite) connected modular matroid with finite rank. The connectivity requirement ensures that every line in a modular space has at least three points.

Note that rank-$3$ modular spaces simplify to projective planes, whose structure is known to exhibit some wild behaviour. However, by Birkoff's Fundamental Theorem on Projective Geometry, modular spaces of rank at least $4$ simplify to projective geometries over division rings.  Nevertheless, we do not need to understand the underlying algebraic structure of modular spaces. 

Henceforth,  {\em matroids} are finite.
A {\em planted matroid} is a pair $(M,Q)$ where $Q$ is a modular space and $M$ is a finite restriction of $Q$. To avoid ambiguity, we will call the flats of $Q$ {\em subspaces} and the closure in $Q$ of a set $Z$ of elements will be called the {\em span} of $Z$, whereas within $M$ we will use the more familiar terms {\em flat} and {\em closure}. We say that two finite sets $V_1$ and $V_2$ in $Q$ are {\em point-wise equivalent} if there is a bijection $\phi:V_1\rightarrow V_2$ such that for each $v\in V_1$ either $v=\phi(v)$, or $\{v,\phi(v)\}$ is a parallel pair, or $v$ and $\phi(v)$ are both loops.

If $(M_1, Q)$ and $(M_2,Q)$ are planted matroids on a common ambient space then we can define 
the intersection and union of $M_1$ and $M_2$ naturally as
\begin{eqnarray*}
 M_1\cap M_2 &:=& Q|(E(M_1)\cap E(M_2)),\mbox{ and }\\
  M_1\cup M_2 &:=& Q|(E(M_1)\cup E(M_2)).
\end{eqnarray*}
For a planted matroid $(M,Q)$ and  $X\subseteq E(Q)$ we denote $M\cup (Q|X)$
by $M\cup X$ and when $X=\{e\}$ we simply write $M+e$. 

Minors behave nicely for planted matroids.
For $e\in E(M)$, the {\em deletion} of $e$ is $(M,Q)\del e := (M\del e,Q)$ and the contraction of $e$ is $(M,Q)\con e:= (M\con e,Q\con e).$
Sometimes we will contract an element $e\in E(Q)\setminus E(M)$, but we do not consider this a minor. We denote $(M\cup\{e\})\con e$ by $M\proj e$, and we call this a 
{\em projection}. For an element $e\in E(M)$, we define the projection, $M\proj e$, to be the matroid obtained from $M\con e$ by adding $e$ back as a loop.

It is worth mentioning, but not important for us, that planted matroids do not behave well under duality, but that is only due to the afore-mentioned issues with rank-$3$ modular spaces. Every rank-$3$ matroid $M$ extends to a modular space, but that need not be the case for the dual.

\section{Decomposition on a separation}

In Section~\ref{modular flats} we considered the composition of matroids using
modular sums; here we consider decomposition.
Let $(M,Q)$  be  a planted matroid. For a matroid $N$ with $E(N)\subseteq E(M)$
we want to determine whether $M$ has $N$ as a minor.
Further suppose that we have a separation $(A,B)$ of order $k$ in $M$ with $E(N)\subseteq A$.
In this section we discuss the issues that arise in decomposing $M$ across the separation $(A,B)$.
For notational reasons it is convenient to assume that $A$ and $B$ are disjoint. If this is not the case,
we will replace the elements of $B$ that are also in $A$ with parallel copies, and leave $A$ unchanged.
Let $M_A=Q|A$ and $M_B=Q|B$.

Let $F$ be the rank-$k$ subspace of $Q$ containing the points in the span of both $A$ and $B$
and let $Z$ be a minimal spanning subset of $F$ such that $Q|Z$ is modular in $M_B\cup Z$.
By Lemma~\ref{modular sets}, each strand $S$ of $M_A$ 
in $M$ spans exactly one element in $Z$ and, by minimality,
each element in $Z$ is spanned by some strand of $M_A$ 
in $M$. Thus $Z$ is finite and uniquely determined up to point-wise equivalence. We call the matroid $Q|Z$ a {\em torso} of  $(A,B)$
in $(M,Q)$. Moreover, since it is uniquely determined up to point-wise equivalence we may also call it {\em the} torso of $(A,B)$ in $(M,Q)$.

Suppose that $T$ is the torso of $(A,B)$ in $(M,Q)$ and
consider decomposing $M$ into $M_A\cup T$ and $M_B\cup T$.
Note that $M\cup T$ is a modular sum of $M_A\cup T$ and $M_B\cup T$. 
It follows from Lemma~\ref{decomp} that, if $M_A\cup T$ does not contain $N$ as a minor,
then neither does $M$. The converse is, however, not as clear-cut.

As noted above,  each element $e$ in the torso is in the span of some strand $S$ for $N$.
We say that $e$ is a {\em base-point} for the strand $S$.
In $M\cup S\cup \{e\}$, the set $S$ is a series-class and when we contract all but one element, say $f$, the pair $\{e,f\}$ will be parallel. 
Thus an element in the torso of $(A,B)$ is any point in the span of $A$ onto which we can contract elements of $M$, without reducing the rank of $M_A$.

One difficulty is that, while the elements in the torso record the points to which we can contract elements, we cannot necessarily contract two prescribed elements down simultaneously. This issue is addressed in graphs by the two-linkage problem, but it is not a problem that we know how to address in matroids. We say that a restriction $T'$ of the torso $T$ is {\em realizable} if there are disjoint sets  $C$ and $Z$ in $B$ such that $r(C\cup A)) = r(C)+r(A)$, and $Z$ is point-wise equivalent with $E(T')$ in $Q\proj C$.

Due to the issues raised by realizability, it is possible that $M_A\cup T$ contains  $N$ as a minor
whereas $M$ does not. 
Nevertheless, we can decompose across  separations of order at most $2$.
\begin{lemma}
Let $T$ be the torso of a separation $(A,B)$ of order at most $2$ in a planted matroid $(M,Q)$
and let $N$ be a matroid with $E(N)\subseteq A$. Then
$N$ is a minor of $M$ if and only if $N$ is a minor of $(Q|A)\cup T$.
\end{lemma}

\begin{proof} Let $M_A=(Q|A)\cup T$.
For the non-trivial direction, suppose that $N$ is a minor of $M_A$ and let $(C,D)$ be a partition
of $E(T)$ such that $N$ is a minor of $M_A\con C\del D$.
We need to prove that $N$ is a minor of $M$.
Since singletons are realizable, we may assume that $C$  contains at least two elements in $T$.
But then $C$ spans $T$ and, hence, $M_A\con C\del D = M_A\con E(T) = M\con B$. So $N$ is a minor of $M$, as required.
\end{proof}

We are, of course, interested in the case that $N\cong U_{2,3}$, for which we can immediately draw the following conclusion.
\begin{lemma}\label{decomp2}
Let $T$ be the torso of a separation $(A,B)$ of order at most $2$ in a planted matroid $(M,Q)$
and let $X\subseteq A$. Then
$X$ is coupled in $M$ if and only if $X$ is coupled in $(Q|A)\cup T$.
\end{lemma}

\section{Rooted matroids}

A {\em rooted matroid} is a triple $(M,Q,V)$ such that $(M,Q)$ is a planted matroid, and $V$ is an independent set in $Q$; we refer to $V$ as the {\em vertex-set} of the rooted matroid. Note that the vertices need not be contained in $M$.
We say that a rooted matroid $(M,Q,V)$ is {\em coupled} if $V$ is coupled in $M\cup V$.
To avoid ugly terminology, we will abbreviate ``coupled rooted matroid" to coupled matroid.

We say that a rooted matroid $(N,Q',V')$ is a {\em minor} of a rooted matroid $(M,Q,V)$ if $(N,Q')$ is a minor of $(M,Q)$ and $V'\subseteq V$.

\subsection*{Bridges}
A {\em bridge} of $(M,Q,V)$ is a rooted matroid $(B,Q,Z)$ such that $E(B)$ is the ground-set of a component of $M\proj V$ and $Z$ is a minimal subset of $V$ such that $\sqcap_Q(Z,E(B)) =\sqcap_Q(V,E(B))$.
Since $V$ is independent, the set $Z$ is uniquely determined; $Z$ is the set of elements in $V$ such that $B\cup Z$ is a component of $B\cup V$. We call $Z$ the set of {\em attachments} of the bridge.  We say that a bridge is {\em trivial} if it is spanned by $V$; if $(B,Q,Z)$ is a trivial bridge, then $|B|=1$.

\subsection*{Torsos of rooted matroids}
The {\em torso} of a rooted matroid $(M,Q,V)$ is the torso of the separation $(V,E(M))$ in $(M\cup V,Q)$.
If $T$ is the torso of $(M,Q,V)$ then we call $(T,Q,V)$ the {\em rooted torso} of $(M,Q,V)$.
 The following lemma is an immediate consequence of Lemma~\ref{modular implies round}.
\begin{lemma}\label{round torso}
If $T$ is the torso of a bridge $(B,Q,X)$ of a rooted matroid $(M,Q,V)$, then $T$ is round.
\end{lemma}

We view a rooted matroid $(M,Q,V)$ as way of building $(M,Q)$ via modular sums.
Let $(B_1,Q,V_1),\ldots,(B_k,Q,V_k)$ be the bridges of $(M,Q,V)$ and let $T$ be the torso.
Now let $T_1,\ldots,T_k$ be the torsos of the bridges $(B_1,Q,V_1),\ldots,(B_k,Q,V_k)$. To avoid ambiguities caused by point-wise equivalence, we take $T_1,\ldots,T_k$ as restrictions of $T$. 
Now the rooted matroid $(M\cup T,Q,V)$ is obtained from  $(T,Q,V)$ by attaching the ``extended bridges'' $(B_1\cup T_1,Q,V_1),\ldots,(B_k\cup T_k,Q,V_k)$ using modular sums; finally, we obtain $(M,Q)$ as a restriction of $(M\cup V\cup T,Q)$.

\subsection*{Frameworks and patchworks}
A {\em framework} is a rooted matroid $(M,Q,V)$ such that each element of $M$ is spanned by a set of at most two elements of $V$. 
Equivalently, $(M,Q,V)$ is a framework if and only if $V$ spans $E(M)$ in $Q$ and $(M,Q,V)$ is coupled.
More generally, a rooted matroid $(M,Q,V)$ is coupled if and only if the rooted torso of $(M,Q,V)$ is a framework.
 
 Frameworks are maintained under taking minors.
 \begin{lemma}\label{minors of frameworks}
 If $(M,Q,V)$ is a framework, $C$ and $D$ are disjoint sets in $M$, and $V'$ is a maximal independent subset of $V$ in $Q\proj C$,
 then $(M\del D\con C,Q\proj C,V')$ is a framework.
 \end{lemma}

 \begin{proof} Inductively it suffices to consider the case that $|C\cup D|=1$, which we leave to the reader.
 \end{proof}
  
We say that a bridge $(B,Q,Z)$ of a rooted matroid $(M,Q,V)$ is a {\em patch}  or a {\em $k$-patch} if $k=|Z|\le 3$ and 
the rooted torso of $(B,Q,Z)$ is a framework.
A rooted matroid $(M,Q,V)$ is a {\em patchwork} if each of its bridges is a patch.  Patchworks are also maintained under taking minors.
 \begin{lemma}\label{minors of patchworks}
 If $(M,Q,V)$ is a patchwork, $C$ and $D$ are disjoint sets in $M$, and $V'$ is a maximal independent subset of $V$ in $Q\proj C$,
 then $(M\del D\con C,Q\proj C,V')$ is a patchwork.
 \end{lemma}

 \begin{proof} Again it suffices to consider the case that $|C\cup D|=1$, which we leave to the reader.
 \end{proof}

Our main result, Theorem~\ref{one}, says that:
{\em if a rooted matroid $(M,Q,X)$ is coupled, then
there is an independent set $V$ in $Q$ containing $X$ such that $(M,Q,V)$ is a patchwork.}

\section{Extending a framework}\label{frameworks}

The results in this section prove our main theorem (Theorem~\ref{one}) in the special case that $M$ is a one-element extension of a framework. 

For a rooted matroid $(M,Q,V)$ and an element $e$ in $Q$ that is spanned by $V$, the {\em attachment-set} of $e$ is the minimal subset of $V$ that spans $e$; since $V$ is independent, this is well-defined.

\subsection*{Skeletons}
Let  $(M,Q,V)$ be a framework. The {\em skeleton} of a framework $(M,Q,V)$ is the graph $G=(V,E')$ defined as follows: $E'$ is the set of non-loop elements of $M$ and the set of ends of an edge $e\in E'$ is the attachment-set of $e$ in the framework. It is common practice to add the loops of $M$ as loop edges of $G$ at arbitrary vertices, but we choose to omit them.

Let $G$ be the skeleton of  $(M,Q,V)$. A cycle $C$ of $G$ is called {\em balanced} if $E(C)$ is a cycle of 
$M$. Let $\cB$ denote the set of all balanced cycles of $G$. The pair $(G,\cB)$ is a {\em biased graph} and $M$ is the {\em frame matroid} determined by $(G,\cB)$; see Zaslavsky~\cite{Zaslavsky}. We assume that the reader is familiar with these notions, but we will review some of the highlights. Let $H$ be a subgraph of $G$ and let $X=E(H)$. We say that $H$ is {\em balanced} if each cycle of $H$ is balanced. If $H$ is balanced then $M|X$ is graphic, but the converse need not hold. The rank of $X$ in $M$ is $|V(H)|-c$, where $c$ is the number of balanced components of $H$.

We note that any element of $V$ that is contained in $M$ will be both a vertex of $G$ and a loop at that vertex. This causes some ambiguity when we write $G-e$, but the meaning will be clear from the context.

As we have already seen, frameworks are preserved under taking minors. 
For any element of $M$, the deletion $(M\del e, Q, V)$ is a framework with skeleton $G-e$.
When $e$ is an edge with distinct ends $a$ and $b$ in $G$, in the contraction $Q\proj e$ the elements $a$ and $b$ will be parallel. So
$(M\con e,Q\proj e,V\setminus \{b\})$ will be a framework whose skeleton is obtained from $G\con e$ by deleting any loops of $M\con e$. Finally consider the case that $e$ is a non-balanced loop at a vertex $a$.
Then $(M\con e,Q\proj e, V\setminus \{a\})$ is a framework; the skeleton $G'$ is obtained from $G$ by first deleting any loops of $M\con e$, then we remove the vertex $a$ from the vertex-set and from the end-set of any edge. (So the non-loop edges that are incident with $a$ in $G$ become loop edges in $G'$.)

For a framework $(M,Q,V)$, let $T$ be the attachment set of an element $e\in E(Q)$ that is spanned by $V$.
For any non-empty subset $X$ of $V$ with $T\cap X\neq\emptyset$, we have $\sqcap_Q(X,(V\setminus X)\cup \{e\})=1$, so, up to point-wise equivalence, there is a non-loop unique element $e'\in Q$ that is spanned by both $X$ and $(T\setminus X)\cup \{e\}$. 
Note that,  $e$ and  $e'$ are parallel in $Q\proj (V\setminus X)$. We call $e'$ the {\em image} of $e$ on $X$.
If $T\subseteq X$, then $e$ and $e'$ will be point-wise equivalent, otherwise, when $X\setminus T\neq\emptyset$, the set $(T\setminus X)\cup\{e,e'\}$ is a circuit.

Now, consider effect, on the attachment set, of contracting an edge $f=xy$ of the skeleton of 
$(M,Q,V)$; here we are allowing $x=y$.
Let $T'$ be the attachment set of $e$ in  $(M\con f, Q\proj f, V\setminus \{y\})$.
If $x=y$, then $T'=T\setminus\{y\}$. Consider the alternative case that $x\neq y$. 
If $T$ does not contain both $x$ and $y$, then 
$T\cup\{e\}$ remains a circuit in $Q\proj f$ and hence $T'=T$. This leaves the more interesting case that $x,y\in T$, which depends
on the location of the image $e'$ of $e$ on the line $\{x,y\}$.
If $f$ is point-wise equivalent with $e'$, then $(T\setminus\{x,y\})\cup\{e\}$ is a circuit in $Q\proj f$, so $T'=T\setminus\{x,y\}$.
Otherwise, when $\{e',f\}$ is independent, $(T\setminus\{x,y\})\cup\{e,e'\}$ remains a circuit in $Q\proj f$ and $\{x,e'\}$ become a parallel pair.
Thus $T'=T\setminus \{y\}$.

Our first result uses the definition of ``$X$-tight" from Section~\ref{graph connectivity}.
\begin{lemma}\label{extend framework 1}
Let $G$ be the skeleton of a framework $(M,Q,V)$, let $T$ be the attachment set of an element $e\in E(Q)\setminus E(M)$ that is spanned by $V$, and let
$X\subseteq V$. If $M\cup X$ is connected, and there is a set $W$ of vertices in $T$  that is $X$-tight, then there is a circuit in $M\cup X\cup \{e\}$ that contains $|W|$ or more elements of $X$.
\end{lemma}

\begin{proof}
Consider a counterexample with $|M|$ minimum; thus no two vertices in $X$ are adjacent.
Let $H$ be a complete graph on $X$.
It is easy to see that $G\cup H$ cannot have a cut-edge, since such an element could be contracted.
 Then, by Lemma~\ref{tight sets 1},
 every edge of $G$ is incident with $W$. Moreover, since $W$ is $X$-tight, we have $V=W\cup X$.
 Furthermore, by minimality, no two vertices of $W$ are adjacent in $G$.
 Finally, by Lemma~\ref{tight sets 2},  each vertex in $W\setminus X$ has degree $2$.
We may assume that $W$ is not contained in $X$ since otherwise the result clearly holds.
Moreover, since $W$ is $X$-tight, there is no circuit in $G$.
Let $x\in X$ be a degree-$1$ vertex of $G$, let $w\in W$ be its neighbour, and let $y\in X\setminus\{x\}$ be the other
neighbour of $w$. Let $f_x$ and $f_y$ be the edges with ends $\{w,x\}$ and $\{w,y\}$ respectively.

Note that $(W\setminus\{w\})\cup \{x\}$ is $X$-tight in $G-w$.  If $x\in T$, then we obtain a contradiction by  contracting $f_y$.
On the other hand, when $x\not\in T$, then we obtain a contradiction by  contracting $f_x$.
\end{proof}

The following result connects tightness with our earlier Sylvester-Gallai-type results. We recall that $K(G,X)$ is the matroid on $V(G)$
such that the rank of a set $Z\subseteq V(G)$ is the connectivity from $Z$ to $X$.
\begin{lemma}\label{extend framework 2}
Let $G$ be the skeleton of a framework $(M,Q,V)$, let $T$ be the attachment set of an element $e\in E(Q)\setminus E(M)$ that is spanned by $V$, and let
$X\subseteq V$. If $G$ is $2$-connected and $K(G,X)|T$ has a three-point plane, then there is a circuit in $M\cup X\cup \{e\}$ that contains three or more elements of $X$.
\end{lemma}

\begin{proof}
Since $G$ is $2$-connected, $K(G,X)$ is simple. Let $P$ be a three-point plane in $K(G,X)|T$. Now let 
{$(A,B)$ be a }$ (P,X)$-separation
in $G$ of order $3$ such that $V(A)$ is the flat of $K(G,X)$ that is spanned by $P$. Note that $V(A\cap B)$ is $X$-tight in $B$
and, since $P$ is a three-point plane in $K(G,X)$, we have $V(A)\cap T=P$. Therefore, we can contract $P$ onto $V(A\cap B)$ and apply
Lemma~\ref{extend framework 1} to get the result.
\end{proof}

The following result is not specifically required, but it does show the strength of the machinery that we have developed.
\begin{lemma}\label{extend framework 3}
Let $G$ be the skeleton of a framework $(M,Q,V)$, let $T$ be the attachment set of an element $e\in E(Q)\setminus E(M)$ that is spanned by $V$, and let
$X\subseteq V$. If $G$ is $3$-connected and $\kappa_{G}(T,X)\ge 4$, then there is a circuit in $M\cup X\cup \{e\}$ that contains three or more elements of $X$.
\end{lemma}

\begin{proof}
Since $G$ is $3$-connected and $\kappa_{G}(T,X)\ge 4$, the gammoid $K(G,X)$ has no triangles and has rank at least $4$.
Since $\kappa_{G}(T,X)\ge 4$, the restriction $K(G,X)|T$ has rank at least $4$. Then, by Motzkin's Lemma, 
$K(G,X)|T$ has a three-point plane. (Motzkin's Lemma, as stated, is for rank-four matroids, but we can reduce to rank four by truncation.)
Now the result follows from Lemma~\ref{extend framework 2}.
\end{proof}

\section{Extending into the torso of a patchwork}

Our main result (Theorem~\ref{one}) is that  every coupled matroid $(M,Q,X)$ extends {to} 
a patchwork. We prove the theorem via induction, resulting in the case that, for some element $e$ in $M$, there is a patchwork $(M\del e,Q,V)$
such that $e$ is spanned by $V$. At that point, it is easy to imagine $(M\del e,Q,V)$ as a framework, allowing
us to apply the results from Section~\ref{frameworks}. However, we need to be careful since there are restrictions of the torso
that are not {realizable}. The results in this section address that problem.

We recall that the rooted torso $(T,Q,V)$ of a coupled matroid $(M,Q,V)$ is a framework.
We refer to the skeleton of $(T,Q,V)$ as the {\em skeleton} of $(M,Q,V)$.

We start by showing that two vertices are adjacent in the skeleton if and only if there is a patch containing the
vertices as attachments; we omit the more obvious direction.
\begin{lemma} \label{cliques in the torso}
If $G$ is the skeleton of a coupled matroid $(M,Q,V)$ and $(B,Q,Z)$ is a bridge of $(M,Q,V)$, then the vertices
in $Z$ are pairwise adjacent in $G$.
\end{lemma}

\begin{proof}
Suppose that there are two distinct non-adjacent vertices $x,y\in Z$ in $G$. Let $T$ be the torso of $(B,Q,Z)$.
By the definition of the attachment set there is an element of $T$ that is not spanned by $Z-x$ and there is an element of $T$ 
that is not spanned by $Z-y$. However, $(T,Q,Z)$ is a framework. Thus we can cover $T$ by two proper flats (those spanned by 
$Z-x$ and $Z-y$), contrary to Lemma~\ref{round torso}.
\end{proof}

We would ideally like to take a minor of a patchwork to get a framework with the same skeleton (up to simplification). This is not always possible, but 
the following result meets our needs. 
\begin{lemma}\label{two out of three}
If $(P,Q,Z)$ is a $3$-patch in a patchwork $(M,Q,V)$, then there is a triangle $T$ in the torso
of $(P,Q,Z)$ such that $Q|(Z\cup T)$ is isomorphic to $M(K_4)$ and each two-element subset
of $T$ is realizable.
\end{lemma}

\begin{proof}
Let $Z=\{x,y,z\}$. Since $P\cup Z$ is connected, there is a circuit $C$ in $P\cup Z$ that contains $y$ and $z$.
Since $Z$ is coupled, $x\not\in C$. Since $P\cup Z$  and $P\proj Z$ are both connected, so is 
$(P\cup\{x\})\proj\{y,z\}$. By Lemma~\ref{strand}, there is a strand $S$ for $(C\setminus\{y,z\})$ in $(P\cup\{x\})\proj\{y,z\}$ that contains $x$.
Let $A=C\cup S$.
Note that 
\begin{itemize}
\item[(i)] $A$ is independent in $Q$,
\item[(ii)] $C$ is a circuit in $Q|A$
\item[(iii)] $Q(A)$ is connected and has corank $2$, and
\item[(iv)] $\sqcap_Q(A\setminus Z,Z)=2$.
\end{itemize}
The torso $T$ of $(A\setminus Z,Q,Z)$ spans a line and by Lemma~\ref{cliques in the torso}, there are at least three elements in $T$.
One of those elements, say $a$, is at the point where the spans of $C\setminus\{y,z\}$ and $\{y,z\}$ meet.
Thus {$\{y,z,a\}$} is a triangle. Since $(A\setminus Z,Q,Z)$ is coupled, every point on $T$ must lie on one of 
the lines spanned by
two of $\{x,y,z\}$, and since $Q|A$ is connected, they cannot all lie on the line spanned by $\{y,z\}$.
Then $T$ must be a triangle and there exist elements $b,c\in T$ such that $\{b,x,y\}$ and $\{c,x,z\}$ are both triangles.
Thus $Q|(T\cup Z)$ is isomorphic to $M(K_4)$. Now  $(Q|A)\con Z$ is connected and has corank $2$, so it
contains a triad, say $A'$. Let $W$ be the matroid obtained from $Q|(A\cup T)$ by contracting 
{$A\setminus A'$}.
Note that $T$ is a modular triangle in $W$ and $A'$ is a triad.  Now by contracting different elements of $A'$,  we see that any two-elements of $T$ is realizable, as required.
\end{proof}

The following result is the upgrade of Lemma~\ref{extend framework 2} to patchworks.

\begin{lemma}\label{extend patchwork 1}
Let $G$ be the skeleton of a patchwork $(M,Q,V)$, let $T$ be the attachment set of an element 
{$e\in E(Q)\setminus E(M)$} that is spanned by $V$, and let
$X\subseteq V$. If $G$ is $2$-connected and $K(G,X)|T$ has a three-point plane, then there is a circuit in $M\cup X\cup \{e\}$ that contains three or more elements of $X$.
\end{lemma}

\begin{proof} Let $W$ be the torso of $(M,Q,V)$. By Lemma~\ref{two out of three}, there is a restriction $W'$ of $W$
such that 
\begin{itemize}
\item[(a)] every pair of adjacent vertices in $G$ remain adjacent in the skeleton of $(W',Q,V)$, and
\item[(b)] every independent subset of $W'$ is realizable in $(M,Q,V)$.
\end{itemize}
By Lemma~\ref{extend framework 2}, there is a circuit $C$ in $X\cup W'$ that contains at least three elements in $X$.
Since $C\setminus X$ is realizable, there is a circuit in $M\cup X\cup e$ that contains at least three elements in $X$.
\end{proof}

We now come to the main result of this section. The statement is a bit convoluted which is in part due to the intended application and 
in part due to the inductive proof. For a graph $G$ we write $\si(G)$ for the simplification of $G$; that is we delete all loops and all but one edge 
from each parallel class.
\begin{lemma}\label{extend patchwork 2}
Let $G$ be the skeleton of a patchwork $(M,Q,V)$, let $T$ be the attachment set of an element 
{$e\in E(Q)\setminus E(M)$} that is spanned by $V$, and let
$X\subseteq V$. If 
\begin{itemize}
\item[(i)] $\kappa_G(T,X)\ge 4$, 
\item[(ii)] the vertices in $X$ are pairwise adjacent in $G$,
\item[(iii)] $G$ is $2$-connected and $\si(G)$ is a subdivision of a $3$-connected graph, and  
\item[(iv)] for each degree-two vertex $v$ in $\si(G)$, there is a $3$-patch $(B,Q,Z)$ of $(M,Q,V)$
such that $v\in Z$ and such that $B$ spans $Z$,
\end{itemize}
then there is a circuit in $M\cup X\cup \{e\}$ that contains three or more elements of $X$.
\end{lemma}

\begin{proof}
Suppose not and consider a counterexample $((M,Q,V),e,T,X)$ with $|V|$ minimum.
Since $G$ is $2$-connected, {$K(G,X)$} is simple. Moreover, since the attachments of 
any $3$-patch are pairwise adjacent, by $(iv)$, there are no two adjacent vertices of degree $2$ in $\si(G)$.
By Motzkin's Lemma, $K(G,X)|T$ has an ordinary plane $P$. By Lemma~\ref{extend patchwork 1},
$P$ is not a three-point plane, so $P$ is the union of a long line $L$ and a point $\{v\}$.
Consider an $(L,X)$-separation $(A,B)$ of order $2$ in $G$ such that $V(A)$ is the closure of $L$ in $K(G,X)$.
Let $x$ and $y$ be the two vertices in $A\cap B$. Every other vertex in $A$ has degree-$2$ in $\si(G)$ and is adjacent to only $x$ and $y$.
Since $L$ is a long line, there are at least three vertices of $T$ in $A$. It follows easily from minimality that every vertex in $A$ is in $T$.

Consider a vertex $z\in V(A)\setminus\{x,y\}$ and let  $(B,Q,Z)$ be a $3$-patch of $(M,Q,V)$
such that $z\in Z$ and such that $B$ spans $Z$. {Note that $Z=\{x,y,z\}$}.
Our goal is to find a strand $S$ in $(B,Q,Z)$, whose base-point  is not in the
the line spanned by $\{x,y\}$, such that the attachment set of $e$ in $(M\con S,Q\proj S,V\del\{z\})$ is $T-\{z\}$.
In doing this, we may lose the edge connecting $x$ and $y$, but that only happens 
in the case that $|L|=3$. In either case $P\setminus \{z\}$ remains an ordinary plane. When $|L|\ge 4$ we have a contradiction to minimality and when $|L|=3$,
we end up with an ordinary plane and can apply Lemma~\ref{extend patchwork 1}.

Let $w$ be the image of $e$ in $Z$ (recall that {\em image} was defined immediately before Lemma~\ref{extend framework 1}).
Now consider the torso $W$ of $(B,Q,Z)$. Note that $W$ is a plane spanning $Z$ and $W$ is round. Since $W$ is round, the set of points
in $W$ outside the span of $\{x,y\}$ still span $W$. Let $I$ be a $3$-element independent subset of $W$ outside the span of $\{x,y\}$.
Since $(B,Q,Z)$ is coupled, the elements of $I$ lie on the two lines spanned by $\{x,z\}$ and $\{y,z\}$. Up to symmetry, we may assume that
two of the three elements in $I$ lie on the line spanned by $\{x,z\}$; one of these two, say $u$, is not on the line spanned by $\{y,w\}$.
Now,  $\{x,y,w\}$ is a triangle in $Q\proj u$. Thus $(T\cup \{e\})\setminus \{z\}$ is a circuit in $Q\proj u$.
So we get the desired minor by taking a strand $S$ in $(B,Q,Z)$ with base point $u$.
\end{proof}

\section{Connectivity in frameworks and patchworks}

Let $(M,Q,V)$ be a patchwork with skeleton $G$ and let $X$ and $Y$ be subsets of $V\cap E(M)$.
Our next goal is to show that the matroid connectivity
$\kappa_{M}(X,Y)$ is equal to the graph connectivity. Recall that our definition of $\kappa_{M}(X,Y)$, in Section~\ref{matroid connectivity}, allows $X$ and $Y$ to intersect. 

\begin{lemma}\label{connectivity through a skeleton}
Let $(M,Q,V)$ be a framework with skeleton $G$ and let $X$ and $Y$ be subsets of $V\cap E(M)$, then the connectivity from $X$ to $Y$ in $M$ is equal to the connectivity from $X$ to $Y$ in $G$.
\end{lemma}

\begin{proof} Let $I$ be the edge set of $k$ disjoint paths from $X$ to $Y$ in $G$, and let $X'$ be the ends in $X$ and $Y'$ be the ends in $Y$.
Then $X'$ and $Y'$ are point-wise equivalent independent sets in the planted matroid $(M\con I,Q\proj I)$, which implies  that $\kappa_{M}(X,Y)\ge k$.
Conversely, if there is a separation $(H_1,H_2)$ of order $k$ in $G$ with $X\subseteq V(H_1)$ and $Y\subseteq V(H_2)$, then no component in the projection $M\proj V(H_1\cap H_2)$ 
contains both an element of $X$ and an element of $Y$, which implies that $\kappa_{M}(X,Y)\le |V(H_1\cap H_2)| = k$.
\end{proof}

A {\em separation} in a rooted matroid $(M,Q,V)$ is a pair $((M_1,Q,V_1),\, (M_2,Q,V_2))$ of rooted matroids such that $M=M_1\cup M_2$, $V=V_1\cup V_2$, and each bridge
of either $(M_1,Q,V_1)$ or $(M_2,Q,V_2)$ is also a bridge in $(M,Q,V)$; the {\em order} of the separation is $|V_1\cap V_2|$.
If $(M,Q,V)$ is a patchwork with skeleton $G$, and $X,Y\subseteq V$, then the {\em connectivity} from $X$ to $Y$ in $(M,Q,V)$ is the minimum order of a separation
$((M_1,Q,V_1),\, (M_2,Q,V_2))$ with $X\subseteq V_1$ and $Y\subseteq V_2$. If $G$ is the skeleton of the patchwork $(M,Q,V)$, then it is easy to see that the connectivity
from $X$ to $Y$ in $G$ is the same as the connectivity from $X$ to $Y$ in $(M,Q,V)$.

\begin{lemma}\label{connectivity through a patchwork}
Let $(M,Q,V)$ be a patchwork and let $X$ and $Y$ be subsets of $V\cap E(M)$, then the connectivity from $X$ to $Y$ in the matroid $M$ is equal to the connectivity 
from $X$ to $Y$ in the patchwork $(M,Q,V)$.
\end{lemma}

\begin{proof} Suppose that the connectivity from $X$ to $Y$ in $(M,Q,V)$ is $k$, and let $G$ be the skeleton of $(M,Q,V)$.
Choose a set of $k$ disjoint paths having as few edges as possible and let $I$ be the edge-set of those paths. By minimality, $I$ contains at most one base point from each patch
of $(M,Q,V)$, therefore the subgraph $G[V,I]$ is realizable, and hence, by Lemma~\ref{connectivity through a skeleton}, we have  $\kappa_{M}(X,Y)\ge k$.
Conversely, there is a separation $((M_1,Q,V_1),\, (M_2,Q,V_2))$ of order $k$ in $(M,Q,V)$ with $X\subseteq V_1$ and $Y\subseteq V_2$, and no non-loop component in the projection $M\proj (V_1\cap V_2)$ 
contains both an element of $X$ and an element of $Y$, which implies that $\kappa_{M}(X,Y)\le |V_1\cap V_2| = k$.
\end{proof}

It is implicit in the proofs of Lemmas~\ref{connectivity through a skeleton} and~\ref{connectivity through a patchwork} that the connectivity
between $X$ and $Y$ in a patchwork is obtained via paths in the skeleton, and that a path $P$ 
{between  distinct} elements $x\in X$ and $Y\in Y$ gives a circuit.
This is made more precise in the following result, whose easy proof is omitted.  
\begin{lemma}\label{paths}
Let $(M,Q,V)$ be a patchwork and let $X$ and $Y$ be subsets of $V\cap E(M)$. If the connectivity from $X$ to $Y$ in $M$ is $k$, then there is a restriction $N$
of $M$ such that $X\cup Y\subseteq E(N)$, $\kappa_N(X,Y)=k$, and each component of $N$ is either a coloop or a circuit containing exactly one element of $X\setminus Y$ and exactly 
one element of $Y\setminus X$.
\end{lemma}

In the rest of this section we consider a minimum order 
{$(X,Y)$-separation} in $M$ relative to the patchwork $(M,Q,V)$.
We show that we can  ``display" the separation in a patchwork but with a possibly different set $V$. As a point of caution, the set $V$ can change in more wild ways than one might imagine upon first reading; for example, as with graphic matroids,  one can create unpleasant examples via Whitney switching. With additional effort it is possible to describe the minimum order $(X,Y)$-separations in $M$, without changing $V$, but, fortunately, we do not need that.
\begin{lemma}\label{making a separation visible}
Let $(M,Q,V)$ be a patchwork, let $X$ and $Y$ be subsets of $V\cap E(M)$, and let $(A,B)$ be a minimum order $(X,Y)$-separation in $M$.
Then there exists $V'\subseteq E(Q)$ and sets $V_1$ and $V_2$ such that 
\begin{itemize}
\item[(a)] $(M,Q,V')$ is a patchwork, and
\item[(b)] $((M|A,Q,V_1),\, (M|B,Q,V_2))$ is a separation of order $\sqcap_M(A,B)$ in $(M,Q,V')$ with $X\subseteq V_1$ and  {$Y\subseteq V_2$}.
\end{itemize}
\end{lemma}

\begin{proof}
Let $Z=X\cap Y$ and let $k=\sqcap_M(A,B) - |Z|$.
By Lemma~\ref{paths}, there exist circuits $C_1,\ldots, C_k$ in $M\del Z$ such that
$r(C_1\cup\ldots\cup C_k\cup Z) = r(C_1) +\cdots+r(C_k)+r(Z)$ and each of the circuits $C_1,\ldots,C_k$
contains exactly one element from each of $X$ and $Y$.
Now let $F$ be the rank-$k+|Z|$ subspace of $Q$ that is spanned by both $A$ and $B$. 
{Note that no non-loop component} of $M\proj F$ contains both an element of $X$ and an element of $Y$.
Therefore, for each $i\in \{1,\ldots,k\}$ there is a non-loop element $v_i\in F$ that is spanned by $C_i$.

Since  $r(C_1\cup\ldots\cup C_k\cup Z) = r(C_1) +\cdots+r(C_k)+r(Z)$, the set $\{v_1,\ldots,v_k\}\cup Z$ is a basis
of $M|F$ and, hence, {each of the circuits} $C_1,\ldots,C_k$ spans exactly one element in $F$.  Therefore, for each $i\in\{1,\ldots,k\}$, the circuit
$C_i$ breaks into two circuits $C_i\cap A$ and $C_i\cap B$ in $M\proj v_i$.

We will show that there exists a basis $V_B$ of the the subspace spanned by $B$ such that
$(B,Q,V_B)$ is a patchwork and, up to point-wise equivalence, $V_B$ contains both
$\{v_1,\ldots,v_k\}\cup Z$ and $Y$. By symmetry, this would also give a basis $V_A$ of the the subspace spanned by $A$ such that
$(A,Q,V_A)$ is a patchwork and, up to point-wise equivalence, $V_A$ contains both
$\{v_1,\ldots,v_k\}\cup Z$ and $Y$. Our result will then follow by taking a basis $V'$ of $V_A\cup V_B$.

Consider the matroid $M_B:=B\cup C_1\cup\ldots\cup C_k$ and the patchwork $(M_B,Q,V)$.
Note that, for each $i\in \{1,\ldots,k\}$, the set  $C_i\setminus A$ is a series class in $M_B$ having $v_i$ as its base point.
Since $V$ is coupled in {$M\cup V$}, this series class contains exactly one element of $V$ and that element, say $x_i$, is in $X$.
When we contract the other elements in the series class, $x_i$ lands in parallel with $v_i$. By Lemma~\ref{minors of patchworks}, we 
get a minor $(B\cup \{x_1,\ldots,x_k\},Q',V_B')$ of the patchwork $(B,Q,V)$ where $V_B'\subseteq V$ and, up to point-wise equivalence,
$V_B'$ contains both $\{v_1,\ldots,v_k\}\cup Z$ and $Y$, as required.
\end{proof}

\section{Extending a patchwork}

Our first result of this section deals with extensions that belong in a patch. 
\begin{lemma}\label{goes into a patch}
Let $(M,Q,V)$ be a patchwork with skeleton $G$, let $X\subseteq V$, and let $e$ be an element of $Q$ in the span of $M\cup V$ such that $(M+e,Q,X)$ is coupled. If $Z\subseteq V$ with $\kappa_G(Z,X)=|Z|\le 3$,  and the connectivity 
from $Z\cup\{e\}$ to $X$ in $M\cup Z\cup X\cup \{e\}$ is $|Z|$, 
then there is a set $V'$ in $Q$, containing $X$, such that $(M+e,Q,V')$ is a patchwork.
\end{lemma}

\begin{proof} Let $k=|Z|$.
Consider a $(Z\cup \{e\}, X)$-separation $(A',B')$ of order $k$ in $M\cup Z\cup X\cup \{e\}$.
Note that {$(A'\del \{e\},B')$} is a minimum order $(Z,X)$-separation in $M\cup Z\cup X$ with $|A'|$ maximum.
This uniquely determines $A'$ and thus $e$ is in the span of $A'$.
By Lemma~\ref{making a separation visible} and possibly changing the set $V$, we may assume that
there is a separation $((A,Q,V_1),(B,Q,V_2))$ of order $k$ in $(M,Q,V)$  such that 
{$Q|A'$ is a restriction 
of $A\cup V_1$ and $Q|B'$ is a restriction of $B\cup V_2$.} 
{Since $k\le 3$}, $(M,Q,V_2)$ is a patchwork. Furthermore, since $e$ is spanned by $A\cup V_1$, each bridge
of $(M+e,Q,V_2)$ has at most three attachments. To complete the proof it remains to show 
{that each bridge} $(B,Q,W)$
of $(M+e,Q,V_2)$ is coupled, moreover it suffices to check the bridge containing $e$ which holds since, by construction, the connectivity from 
$W$ to $X$ in the skeleton is $k=|W|$, so we can take a minor of  $(B,Q,V_2)$ such that $W$ is point-wise equivalent to a subset of $X$.
\end{proof}

Roughly the next result says that, if we ``block" a patch when extending a patchwork, and the extension remains coupled, then we can grow the vertex-set of our patchwork such that $e$ looks like an edge from the point of view of the patch.
\begin{lemma}\label{blocking a patch}
Let $(M,Q,V)$ be a patchwork with skeleton $G$, let $X\subseteq V$, and let $e$ be an element of $Q$ in the span of $M\cup V$ such that $(M+e,Q,X)$ is coupled.
Now, for $k\in\{2,3\}$, let $(P,Q,Z)$ be a patch, with $\kappa_G(Z,X)=|Z|=k$, such that
$e$ is not in the span of $(E(M)\cup V)\setminus E(P)$ and $\kappa_{M\cup Z}(P\cup Z\cup \{e\},X)=k+1$. Then there exists $v$ in the span of $P\cup Z$
and $w$ in the span of $(E(M)\cup V)\setminus E(P)$ such that $e$ is on the line spanned by $v$ and $w$ in $Q$ and such that 
$(P,Q,Z\cup\{v\})$ is coupled.
\end{lemma}

\begin{proof}
By Tutte's linking theorem, there is a minor $(P\cup Z\cup\{e\},Q',X)$ of $(M\cup Z\cup\{e\},Q,V)$ such that the local connectivity between
$P \cup Z\cup\{e\}$ and $X$ is $k+1$ and the local connectivity between $Z$ and $X$ is $k$ in $Q'$.
Since $Z$ and $X$ are both contained in $V$ and $V$ is coupled in $(M,Q,V)$, the set $Z$ is point-wise equivalent with some
subset $X'$ of $X$ in $Q'$. Let $D$ be a basis of $P\cup X$ that contains $X$. There is an element $v$ in $Q'$, in the span of $P\cup Z$
such that $(D\setminus X')\cup \{v\}$ spans $\{e\}$; moreover this point is unique, up to equivalence. 
Note that the span of $P\cup Z$ in $Q'$ contains the span of $P\cup Z$ in $Q$, so we may assume that the element $v$ 
is in the span of $P\cup Z$ in $Q$. 

Now, up to equivalence, there is a unique point $w$ in $Q'$ that is spanned by $X\setminus X'$ such that $e$ is on the line spanned by 
$v$ and $w$ in $Q'$. Note that $Q'$ contains the line spanned by $v$ and $e$ in $Q$ so we can choose $w$ on that line. Since $Q'$ is 
obtained from $Q$ by projecting out elements in $E(M)\setminus E(P)$, the point $w$ is spanned by $(E(M)\cup Z)\setminus E(P)$ in $Q$.
Since $P\cup Z$ is connected, we can contract $v$ into the span of $X'$. However, $(P\cup\{e\},Q',X)$ is coupled, so 
$w'$ must be point-wise equivalent with an element of $X\setminus X'$ and, however, we contract $v$ into the span of $X'$ it must land in parallel with an element of $X'$. Thus $(P,Q',Z\cup\{v\})$ is coupled, and hence, so is $(P,Q,Z\cup\{v\})$.
\end{proof}

\section{Resilient patchworks}

In the proof of the main theorem, it is convenient to work with a patchwork $(M,Q,V)$ such that
\begin{itemize} 
\item[(1)] for each patch $(P,Q,Z)$ with $|Z|=3$, the matroid $P\cup Z$ is vertically $3$-connected, and
\item[(2)] the skeleton of $(M,Q,V)$ is $3$-connected.
\end{itemize}
We say that the patchwork is {\em resilient} if it satisfies these two properties.

\begin{lemma}\label{connectivity in a 3patch}
Let $(P,Q,Z)$ be a $3$-patch in a patchwork $(M,Q,V)$.
If $M\cup V$ is vertically $3$-connected, but $P\cup Z$ is not vertically $3$-connected, then there is 
{an element $v$ in $Q$} that
is spanned by $P\cup Z$, but not by $V$, such that  $(M,Q,V\cup\{v\})$ is a patchwork.
\end{lemma}

\begin{proof}
Consider a vertical separation $(A,B)$ of order $1$ in $P\cup Z$. 
Since $M\cup V$ is vertically $3$-connected, there is an element $z_1\in E(A)\cap Z$ that is not spanned by $B$
and an element $z_2\in E(B)\cap Z$ that is not spanned by $A$. Let $v$ be a non-loop element in $Q$ that 
is spanned by both $A$ and $B$. Since $(P,Q,Z)$ is a patch, $P\proj Z$ is connected, and hence $v$ is not spanned by $Z$.
It is now easy to see that $(M,Q,V\cup\{v\})$ is a patchwork, as required.
\end{proof}

{Therefore, if we are given a patchwork $(M,Q,V')$ such that
$M\cup V'$ is vertically $3$-connected,} we can obtain (1) by adding vertices.
If the skeleton of {the resulting patchwork $(M,Q,V)$} is not $3$-connected, then we will get 
a separation $((A,Q,V_1),(B,Q,V_2))$ of order at most $2$ in $(M,Q,V)$ with $V_1\setminus V_2$ and 
$V_2\setminus V_1$ both non-empty. Note that, $(A,Q,V_1)$ and $(B,Q,V_2)$ are both
patchworks, and, moreover, if $(M,Q,V)$ satisfies $(1)$ then so do both $(A,Q,V_1)$ and $(B,Q,V_2)$.

We have proved the following result.
\begin{lemma}\label{resilient}
Let $(M,Q,X)$ be a rooted matroid such that $M\cup X$ is vertically $3$-connected and
each pair of distinct elements in $X$ span a long-line in $M\cup X$. If there is a patchwork
whose vertex-set contains $X$, then there is a resilient patchwork whose vertex-set contains $X$.
\end{lemma}

\section{The main theorem}

Our main result is:
\begin{theorem}\label{one}
If $(M,Q,X)$ is a coupled matroid, then there is a patchwork $(M,Q,V)$ with $X\subseteq V$.
\end{theorem}  

\begin{proof}
We consider a counterexample $(M,Q,X)$  such that
$r(M)-r(X)$ is minimum, and, subject to that, the number of elements of $M$ that are not spanned by $X$ is minimum.
We may further assume that: $M\cup V$ is simple, $X\cap E(M)$ is empty, and each pair of distinct elements in 
$X$ span a long line in $M\cup X$.
Note that $|X|\ge 4$, since otherwise $(M,Q,X)$ is, itself, a patchwork.
By Lemma~\ref{decomp2},
$M\cup X$ is vertically $4$-connected.

Let $e$ be an element of $M$ that is not in the span of $X$. 
By minimality, there is an independent set $V$ in $Q$, containing $X$, such that
$(M\del e,Q,V)$ is a patchwork. We choose $V$ such that
$(M\del e,Q,V)$ is resilient (see {Lemma~\ref{resilient}}) and, subject to that, with $|V|$ is as large as possible.
We extend $V$ to a basis $I$ of $(M\del e)\cup V$ and let $W$ be a minimal subset of $I$ that spans $e$.

By Lemma~\ref{goes into a patch}, we have
\begin{claim}
For any set $Z\subseteq V$ with $|Z|\le 3$,  the connectivity 
from $Z\cup\{e\}$ to $X$ in $M\cup Z\cup X\cup \{e\}$ is $|Z|+1$.
\end{claim}

We will call a patch $(B,Q,Z)$ of $(M\del e,Q,V)$ {\em blocked} if {$E(B)\cap W$} is non-empty.

Our first non-trivial claim is:
\begin{claim}\label{complete the frame}
If $(B,Q,Z)$ is a blocked patch of $(M\del e,Q,V)$, then $|Z|=2$
and there exists  an element $v$ in $Q$ that is spanned by $B\cup Z$ such that
$(M,Q,V\cup\{v\})$ is a patchwork and such that $e$ is spanned by $(W\setminus E(B))\cup \{v\}$.
\end{claim}

\begin{proof}[Proof of claim.]
By Lemma~\ref{goes into a patch}, we may assume that the connectivity 
from $Z\cup\{e\}$ to $X$ in $M\cup Z\cup X\cup \{e\}$ is $|Z|+1$. Then, by Lemma~\ref{blocking a patch},
for each patch $(B,Q,Z)$ of $(M\del e,Q,V)$ there is an element $v\in E(Q)$ that is spanned by $B\cup\{Z\}$ such that
$(B,Q,Z\cup\{v\})$ is coupled and $e$ is spanned by $(W\setminus E(B))\cup \{v\}$.
If $|Z|=2$, then $(M\del e,Q,V\cup\{v\})$ is a patchwork, as required. Therefore, we may assume that $|Z|=3$.

Let $B^+$ be obtained from $B$ by adding elements so that the lines spanned by any two elements in $Z\cup\{v\}$ in
$B\cup Z\cup\{v\}$ are long.
By minimality, there is a patchwork $(B,Q,V')$ with $Z\cup\{v\}\subseteq V'$ and we can choose $V'$ 
such that $(B,Q,V')$ is resilient.
Now $(M\cup B',Q,V\cup V')$ is a resilient patchwork, but, by our choice of $V$, the patchwork $(M,Q,V\cup V')$
is not resilient. Now $(M\cup B',Q,V\cup V')$ and $(M,Q,V\cup V')$ have the same non-trivial patches,
so the issue is that the skeleton $H$ of $(M,Q,V\cup V')$ is not $3$-connected.
However, since $(M,Q,V)$ is resilient, $B\cup Z$ is vertically $3$-connected, so the skeleton of
$(B,Q,V')$ is $2$-connected. Moreover, there are three disjoint paths from $Z$ to $X$ in $H$. So there are three
internally-disjoint paths connecting any two vertices of $Z$ in $H$. Consider a separation $(H_1,H_2)$ of order $2$
in $H$ with $X\subseteq V(H_1)$ and $V(H_2)\setminus V(H_1)\neq \emptyset$, and, subject to that take $|V(H_2)|$ maximum.
Note that $Z\subseteq V(H_1)$ and $v\in V(H_2)\setminus V(H_1)$. Moreover $V(H_2)$ is the closure
of the neighbour set of $v$ in $H$ in the matroid $K(H,X)$. Now $(M,Q,V(H_1))$ is a resilient patchwork with
$V\subseteq V'$, however, $V\neq V'$ since $v$ is in the span of a $3$-patch $(B,Q,Z)$ in $(M,Q,B)$,
whereas in $(M,Q,V(H_1))$ we can disconnect $v$ from $X$ by projecting out only two vertices, namely, $V(H_1\cap H_2)$.
\end{proof}

Applying Claim~\ref{complete the frame} to each blocked patch of $(M\del e,Q,V)$ we obtain a set $V^+$ such that 
\begin{itemize}
\item[(a)] $(M\del e,Q,V^+)$ is a patchwork with $V\subseteq V^+$,
\item[(b)] $V^+$ spans $e$, 
\item[(c)] for each blocked patch $(B,Q,Z)$ of $(M\del e,Q,V)$,
there is exactly one element of $V^+\setminus V$ that is spanned by $B\cup Z$, and
\item[(d)] for each $v\in V^+\setminus V$, there is a blocked patch $(B,Q,Z)$ of $(M\del e,Q,V)$
such that $v$ is spanned by $B\cup Z$.
\end{itemize}

Since $M\cup X$ is vertically $4$-connected, each nontrivial patch $(B,Q,Z)$ {of $(M\del e,Q,V^+)$}
has three attachments and $B$ spans $Z$.
Let $G^+$ denote the skeleton of $(M\del e,Q,V^+)$ and let $T$ be the attachment set of $e$.
Note that
\begin{itemize}
\item[(i)] each vertex $v\in V^+\setminus V$ is contained in $T$ and has exactly two neighbours in $G^+$,
\item[(ii)] any two distinct vertices in $V$ are adjacent in $G$ if and only if they are either
adjacent in $H$ or they have a common neighbour in $V^+\setminus V$ in $H$.
\end{itemize}
This implies that $\si(G^+)$ is a subdivision of a $3$-connected graph. Now it follows from 
Lemma~\ref{extend patchwork 2}, that $\kappa_{G^+}(T,X)\le 3$. But then we get a contradiction by applying
Lemma~\ref{goes into a patch}.
 \end{proof}

\section*{Statement of AI use}
Microsoft Copilot was used to proofread the final draft and to complete the bibliographic citations.

\section*{Acknowledgements}
The matroid minors papers are long overdue, and we are deeply grateful to Neil Robertson for giving us a gentle prod
to resume writing. That prod was accompanied by a very generous personal financial donation from Neil to the University of Waterloo to assist with the writing.

\end{document}